%% file: ray.tex
\documentclass[12pt,fleqn,leqno]{article}
\title{Ray class fields of global function fields with many rational places}
\author{Roland Auer}
\date{}

\usepackage{latexsym}
\usepackage{delarray}
\usepackage{amssymb}
\usepackage[fleqn,leqno]{amsmath}
\usepackage{amsthm}
\swapnumbers

\newlength\tablesep \tablesep0.5cm
\newlength\stateindent \stateindent0.5em

\newcounter {state}[section]
\newcounter {substate}[state]

\renewcommand {\thesubstate} {(\alph{substate})}

\newcommand {\comment}[1] {}
\newcommand {\df}[1] {{\bfseries #1}}
\newcommand {\eg} {\mbox{e.\ g}.\ }
\newcommand {\Eg} {\mbox{E.\ g}.\ }
\newcommand {\ie} {\mbox{i.\ e}.\ }
\newcommand {\wrt} {\mbox{w.\ r.\ t}.\ }
\newcommand {\wl} {\mbox{w.\ l.\ o.\ g}.\ }
\newcommand {\cf} {cf.\ }
\newcommand {\resp} {resp.\ }

\newcommand {\C} {\mathcal C}
\newcommand {\Cl} {\mathcal C\!\ell}
\newcommand {\D} {\mathcal D}
\newcommand {\F} {\mathbb F}
\newcommand {\I} {\mathcal I}
\newcommand {\M} {\mathcal M}
\newcommand {\N} {\mathbb N}
\newcommand {\Norm} {\mathcal N}
\renewcommand {\O} {\mathcal O}
\newcommand {\Q} {\mathbb Q}
\newcommand {\R} {\mathbb R}
\newcommand {\Z} {\mathbb Z}

\renewcommand {\d}{\mathfrak d}
\newcommand {\eps} {\varepsilon}
\newcommand {\f}{\mathfrak f}
\newcommand {\la} {\lambda}
\newcommand {\m}{\mathfrak m}
\renewcommand{\o}{\mathfrak o}
\newcommand {\om} {\omega}
\newcommand {\p} {\mathfrak p}
\newcommand {\ph} {\varphi}
\newcommand {\q} {\mathfrak q}
\renewcommand {\th} {\vartheta}

\newcommand {\ab} {{\rm ab}}

\newcommand {\brackets} {\makebox[1.2em]{$[\;\;]$}}
\newcommand {\floor}[1] {\lfloor #1 \rfloor} 
\newcommand {\into} {\hookrightarrow}
\newcommand {\noqed}{\renewcommand{\qed}{}}
\newcommand {\noresym}[2] {(\;\;,#1|#2)}

\newcommand {\onto} {\twoheadrightarrow}
\newcommand {\on}[2]{\genfrac{}{}{0pt}{}{#1}{#2}}

\newcommand {\rank} {{\rm rank\,}}
\newcommand {\round}[2] {\genfrac{\lceil}{\rceil}{0pt}1{#1}{#2}}

\newcommand {\thru}[3] {#1_{#2},\ldots,#1_{#3}}
\newcommand {\set}[2] {\{#1,\ldots,#2\}}
\newcommand {\sothat} {\,|\,}

\newcommand {\ti}[1] {\tilde{#1}} 

\DeclareMathOperator{\reg}{reg}

\newenvironment{map}[1] 
  {\[#1:\begin{array}{rcl}} 
  {\end{array}\]
  \\[-0.5\baselineskip]  
}

\newenvironment{substate}
{
  \begin{list}{\bf\thesubstate}
    {\usecounter{substate} 
     \itemindent0em
     \settowidth\labelwidth{\bf(b)} \labelsep0.5em  
     \leftmargin\labelwidth \addtolength\leftmargin\labelsep
     \topsep0.5ex 
     \itemsep0ex}
}
{\end{list}}

\newenvironment {subproof}
{
  \begin{list}{\bf\thesubstate}
    {\usecounter{substate} 
     \leftmargin0em 
     \settowidth\labelwidth{\bf(b)} \labelsep0.5em  
     \itemindent\labelwidth \addtolength\itemindent\labelsep
     \topsep0.5ex 
     \itemsep0ex}
}
{\end{list}}

\newtheorem{thm}{Theorem}[section]
\newtheorem{prop}[thm]{Proposition}
\newtheorem{lem}[thm]{Lemma}
\newtheorem{cor}[thm]{Corollary}

\newtheorem{conlem}[thm]{Conductor Lemma}
\newtheorem{hayesex}[thm]{Example (Ray Class Fields \`a la Hayes)} 
\newtheorem{oesthm}[thm]{Theorem (Oesterl\'e)} 
\newtheorem*{diffor}{Different Formula}

\theoremstyle{definition}
\newtheorem{ex}[thm]{Example}

\begin{document}
  \maketitle
  \begin{abstract}
    A general type of ray class fields of global function fields is 
    investigated. 
    The systematic computation of their genera leads to new examples of 
    curves over finite fields with comparatively many rational points. 
  \end{abstract}
  { \renewcommand\thefootnote{}
  }
  \setcounter{section}{-1}
  \input{intro.tex}
  \input{rayfield.tex}
  \input{bounds.tex}
  \input{many.tex}

  \input{rami.tex}

\input{ref.tex}
  \input{address.tex}
\end{document}

%% file: intro.tex
\section{Introduction}

Algebraic curves over finite fields with many rational points have been of 
increasing interest in the last two decades. 
The question of explicitly determining the maximal number of points on a 
curve of given genus was inititiated and in some special cases solved by 
Serre \cite{S1,S2,S3} arround 1982. 
Since then it has been tried to attack the problem by means of algebraic 
geometry as well as field arithmetic. 
Constructions by explicit equations have been carried out by van der Geer and van der Vlugt \cite{GV0,GV3}.
The present paper, which makes use of class field theory, has its 
immediate predecessors in work by Lauter \cite{Lt1,Lt2,Lt3} and 
Niederreiter and Xing \cite{NX3,NX4,NX5,NX9,NX10,XN}. 
The numerical results obtained improve several entries of the tables given in \cite{GV7}, \cite{Ni} and \cite{NX11}. 

As we are looking from the field theoretic point of view,
to an algebraic curve $X$ (smooth, projective, absolutely irreducible) 
defined over a finite field $\F_q$ we associate its field $K=\F_q(X)$ of 
algebraic functions, a global function field with full constant field 
$\F_q$. 
Its genus is that of $X$, and coverings of $X$ correspond to field 
extensions of $K$, the degree of the covering being the degree of the 
extension.

A place of $K$, by which we mean the maximal ideal $\p$ in some discrete 
valuation ring of $K$, with (residue field) degree $d = \deg\p$, 
corresponds to (a Galois conjugacy class of) $d$ points on $X(\F_{q^d})$, 
and each point on $X$ having $\F_{q^d}$ as its minimal field of definition over $\F_q$ lies in such a conjugacy class. 
In particular the rational places, \ie the places of degree 1, of $K|\F_q$ 
are in 1--1 correspondence with the $\F_q$-rational points on $X$.
The (normalized) discrete valuation associated to a place $\p$ of $K$ will 
be denoted by $v_\p$.

Our access to the problem of constructing curves with many rational points will be by ray class fields:
To each non-empty set $S$ of places and each effective divisor $\m$ of 
$K$ with support disjoint from $S$ we associate a global function field 
$K_S^\m$, the largest abelian extension of $K$ with conductor $\leq\m$ in which every place of $S$ splits completely. 
These ray class fields have been used by Perret~\cite{Pe} before, but because of some misleading statements in his paper, we shall give thorough proof for their existence and derive some of their basic properties, a task which is performed in Section~1. 

Section~2 discusses upper bounds for $S$-class numbers and for the number of rational places of a global function field. 
In Section~3 we prove the existence of curves with many rational points by explicitly computing the needed invariants of certain $K_S^\m$. 
At this point I want to thank the referee for drawing my attention to 
recent articles by Niederreiter and Xing, which have appeared after submission of this paper and could be added to the references during the revision.

\comment{
The $K_S^\m$ provide a rather fine classification of the finite abelian 
extensions of $K$ and are canonical and quite beautiful objects in themselves. 
Two interesting and to a large extent still unsolved questions are: 
What are their degrees, genera and numbers of rational places? 
Is there a way to find for each of them an (explicit set of) algebraic function(s) generating it over $K$? 
The present paper attempts an answer to the first question in a more
general situation than this has been done before. 
As for the second problem there is some hope to find a solution in terms of composites of Kummer and Artin-Schreier extensions as outlined in~\cite{GV3}. 
}

%% file: rayfield.tex
\section{Ray Class Fields}

Let $K$ be a global function field with full constant field $\F_q$. 
Denote by $\I = \I_K$ and $\C = \C_K = \I/K^*$  the id\`ele group and the 
id\`ele class group of $K$, respectively. (The diagonal embedding of 
$K^*$ into $\I$ is considered as inclusion.) 

Given a place $\p$ of $K$, let $K_\p$ be the completion of $K$ at $\p$, and
$\bar\p$ the extension of $\p$ to $K_\p$ (\ie the topological 
closure of $\p$ in $K_\p$).
Furthermore we have the $n$-th one-unit group $U_\p^{(n)} = 1 + 
\bar\p^n$, where we admit $n=0$, too, in which case we mean the usual
unit group $U_\p = \{ u \in K_\p^* \sothat v_{\bar\p}(u)=0 \}$.
The canonical embedding of $K_\p^*$ into $\I$ will be expressed by the 
symbol $\brackets_\p$, \ie $[z]_\p$ for $z\in K_\p^*$ is the id\`ele 
$(\ldots,1,z,1,\ldots)$ having $z$ at its $\p$-th position and $1$ 
elsewhere.

Let $S$ be a set of places of $K$. 
By an \df{$S$-cycle (of $K$)} we mean an effective divisor of $K$ with 
support disjoint from $S$. 
Let $\m = \sum m_\p \p$ be such an $S$-cycle. 
Then, following Perret~\cite[pp.~305f]{Pe}, we define the \df{$S$-congruence subgroups mod $\m$},
\[
  \I_S^\m := \left( \prod_{\p\in S}K_\p^* 
             \times \prod_{\p\notin S}U_\p^{(m_\p)} \right) \cap \I 
  \quad\mbox{and}\quad \C_S^\m := K^*\I_S^\m / K^*,
\]
of $\I$ and $\C$, respectively. (Here intersection with $\I$ is only 
necessary, of course, when $S$ is infinite.)
We also put $\I^\m := \I^\m_{S(\m)}$ where $S(\m)$ denotes the set of all 
places of $K$ not occurring in the support of $\m$.

For the purposes of class field theory, $\I$ carries the restricted 
product topology of the $K_\p^*$ with respect to the $U_\p$ (\cf 
\cite[p.\ 68]{CF}), which amounts to saying that $\I$ is a topological 
group and the congruence subgroups $\I^\m_\emptyset$, where $\m$ runs 
through all effective divisors of $K$, form a basis of open 
neighbourhoods of $1$.
By definition the $\I_S^\m$ are open (and thereby closed) subgroups of 
$\I$, and since $\C$ is endowed with the quotient topology induced from 
$\I$, the $\C_S^\m$ are open (and closed), too.

For any place $\p$ of $K$ the kernel of the map $K_\p^* \to
\I/K^*\I_\emptyset^\m$ induced by the embedding $\brackets_\p$ is
contained in $\F_q^*U_\p^{(m_\p)} \subseteq U_\p$. But $K_\p^*/U_\p
\simeq \Z$, hence the index of $\C_\emptyset^\m$ in $\C$ turns out to be infinite, contrary to the respective assertions in \cite[p.\ 302]{Pe}, which are therefore wrong. 
Due to this fact we assume $S$ non-empty from now on. 

Denote by $\O_S$ the ring of all functions in $K$ with poles only in 
$S$, which is known to be a Dedekind domain with quotient field $K$. 
As is easily verified by means of (weak) approximation, the \df{$S$-ray 
class group \mbox{mod $\m$}}, $\C/\C_S^\m\simeq \I/K^*\I_S^\m$, is 
naturally isomorphic to the \df{$S$-ideal class group \mbox{mod $\m$}}, 
$\Cl^\m(\O_S)$, defined as the quotient of the group of fractional ideals 
of $\O_S$ prime to $\m$ by its subgroup of principal ideals $x\O_S$ with 
$x \in K^* \cap \I^\m$. 
Note that for $\m=\o$, the zero divisor of $K$, we regain the usual ideal 
class group $\Cl(\O_S)$, whose cardinality, the \df{$S$-class number (of 
$K$)}, is denoted $h_S$. 

\begin{prop}\label{idelprop}
  The sequence
  \[
    \O_S^* \to \I_S^\o/\I_S^\m \to \I/K^*\I_S^\m \to \I/K^*\I_S^\o \to 1
  \]
  is exact, and $(\C:\C_S^\m)$ is finite.
\end{prop}

\begin{proof}
  We have $\O_S^* = K^* \cap \I_S^\o$, implying $\O_S^*\I_S^\m / \I_S^\m 
  = K^*\I_S^\m\cap\I_S^\o / \I_S^\m$, which is the exactness at 
  $\I_S^\o/\I_S^\m$. The exactness at the other groups is even simpler. 
  As for the second claim note that both $\I_S^\o/\I_S^\m 
  \simeq \prod_\p U_\p/U_\p^{(m_\p)}$ and $\I/K^*\I_S^\o \simeq \Cl(\O_S)$
  are finite.
\end{proof}

Let us briefly recall the essential facts of (global) class field theory 
(see \cite[p.\ 172]{CF} and \cite[p.\ 408]{Nk}): 
For each finite Galois extension $L|K$ there is a surjective homomorphism 
\[ 
  \noresym LK : \C_K \to G(L|K)^\ab
\]
from the id\`ele class group of $K$ onto the abelianization of the Galois 
group of $L|K$, called the \df{norm residue symbol of $L|K$}, which at the 
same time we think of as defined on $\I_K$. 
Its kernel is the \df{norm group} $\Norm_L := N_{L|K} \C_L = 
K^*N_{L|K}\I_L/K^*$, where $N_{L|K}$ denotes the norm from $L$ to $K$. 
Furthermore, the map $L \mapsto \Norm_L$ is an inclusion-reversing 
1--1 correspondence between the finite abelian extensions (within some 
fixed algebraic closure $\bar K$) of $K$ and the closed subgroups of finite 
index in $\C_K$. 

Accordingly we can define the \df{$S$-ray class field \mbox{mod $\m$}}, 
denoted $K_S^\m$, as being the unique (in $\bar K$) abelian extension 
$L|K$ satisfying $\Norm_L = \C_S^\m$. 
Then $K_S^\o$ is just the Hilbert class field of $\O_S$ in the sense of 
Rosen \cite{Ro2} (with the slight generalization that we allow $S$ to be 
infinite). 
Moreover, in case $S$ consists of exactly one place, an explicit 
construction of $K_S^\m$ via rank 1 Drinfel'd modules has been carried out 
by Hayes \cite{Ha}. 

Next we want to investigate some extremality properties of the $S$-ray class fields. 
To this end we introduce the \df{conductor} of a finite abelian extension
$L|K$ as the effective divisor 
\[
  \f(L|K) := \sum_\p f(L,\p) \p 
\]
of $K$, with the \df{conductor exponent $f(L,\p)$ of $\p$ in $L$} being the least integer $n\geq0$ such that the \df{$n$-th upper ramification group $G^n(L,\p)$ of $\p$ in} $L$ gets trivial. 
For a definition of the upper (index) ramification groups we refer the 
reader to \cite[pp.\ 33ff]{CF}, \cite[pp.\ 186ff]{Nk} or the appendix of 
this paper. 
By introducing the conductor in this way we avoid dealing with infinite class field extensions as encountered in Perret's paper~\cite{Pe}. 
For the following lemma, proposition and theorem \cf Section~I of that paper. 
Because $G^0(L,\p)$ is the inertia group of $\p$ in $L$, exactly the places occurring in $\f(L|K)$ are ramified in $L$. 
\begin{conlem}\label{conlem}
  Let $L|K$ be finite abelian and let $S$ be any non-empty set of places    of $K$ splitting completely in $L$. 
  Then $\f(L|K)$ is the smallest $S$-cycle $\m$ such that $L$ is 
  contained in $K_S^\m$ (\ie $\Norm_L$ contains $\C_S^\m$).
\end{conlem}

\begin{proof}
  For any place $\p$ of $K$ the composite map $K_\p^* \into \I \onto   
  G(L|K)$ of the embedding $\brackets_\p$ with the norm residue symbol 
  $\noresym LK$ has the decomposition group $G^{-1}(L,\p)$ as its image and 
  throws $U_\p^{(n)}$ onto $G^n(L,\p)$ for all $n\in\N_0$ 
  (\cf \cite[p.\ 155]{CF}). 
  Hence we have 
  \[
    [U_\p^{(n)}]_\p \subseteq K^*N_{L|K}\I_L \iff n \geq f(L,\p)
  \]
  and even $[K_\p^*]_\p \subseteq K^*N_{L|K}\I_L$ when $\p\in S$. 
  
  Now let $\m = \sum m_\p \p$ be an $S$-cycle of $K$. $\I_S^\m$ is 
  topologically generated by its subgroups $[K_\p^*]_\p$, $\p\in S$, and 
  $[U_\p^{(m_\p)}]_\p$, $\p\notin S$. But $K^*N_{L|K}\I_L$ is closed, hence
  the above facts may be summarized as 
  \[
    \I_S^\m \subseteq K^*N_{L|K}\I_L \iff \m\geq\f(L|K),
  \]
  which is what we had to show.
\end{proof}

For a non-empty set $S$ of places of $K$ denote by $\deg S$ the greatest 
common divisor of the degrees of its elements. 
\begin{prop}\label{rayprop}
  Let $S$ be a non-empty set of places and $\m$ an $S$-cycle of $K$. 
  Then $K_S^\m$ is the largest abelian extension $L$ of $K$ such that 
  $\f(L|K)\leq\m$ and that every place in $S$ splits completely in $L$.
  Moreover, $\F_{q^d}$ with $d := \deg S$ is the full constant field of 
  $K_S^\m$.
\end{prop}

\begin{proof}
  The maximality of $K_S^\m$ \wrt the imposed conditions follows directly 
  from the lemma (and its proof). 
  As for the second assertion, consider the constant field extension 
  $K_d = \F_{q^d}K$ of degree $d$ over $K$. 
  It is unramified and such that all the places in $S$ are completely 
  decomposed. 
  Hence, by the maximality just proven, $\F_{q^d} \subseteq K_d \subseteq 
  K_S^\o \subseteq K_S^\m$. 
  On the other hand, in order that $\F_{q^{d'}}$ for some $d'\in\N$ is 
  contained in $K_S^\m$, $d'$ must divide the degree of every place 
  splitting completely in $K_S^\m$. 
\end{proof}

For a finite abelian extension $L|K$, a place $\p$ of $K$ and $n\geq-1$ we 
denote by $L^n(\p)$ the subfield of $L$ fixed by $G^n(L,\p)$ and call it 
the \df{$n$-th upper ramification field of $\p$ in $L$}.
Given a set $S$ of places of $K$, any effective divisor $\m$ of $K$ can 
always be made an $S$-cycle, denoted $\m \setminus S$, by removing all the 
places in $S$ from its support. 
With this notation we can write down the remarkable fact that the 
ramification fields of the $S$-ray class field mod $\m$ are the $S$-ray 
class fields mod $\m'$ for certain $\m'\leq\m$.
\begin{thm}\label{raythm}
  Let $S$ be a non-empty set of places of $K$, $\m = \sum_\p m_\p\p$ an 
  $S$-cycle and $L$ an intermediate field of $K_S^\m|K$.
  \begin{substate}
  \item For any place $\p$ of $K$ and any $n\in\N_0$ we have 
    \[
      L^n(\p) = L \cap K_S^{\m\setminus\{\p\}+n\p} \mbox{ and } 
      L^{-1}(\p) = L \cap K_{S\cup\{\p\}}^{\m\setminus\{\p\}}
                 = L \cap K_{\{\p\}}^{\m\setminus\{\p\}}.
    \]
  \item The discriminant $\d(L|K)$ of $L|K$ satisfies
    \[
      \d(L|K) = [L:K]\m - \sum_\p \left(
        \sum_{n=0}^{m_\p-1}[L^n(\p):K] 
      \right) \p.
    \]
  \end{substate}
\end{thm}

\begin{proof}
  \begin{subproof}
  \item Consider an arbitrary intermediate field $L'$ of $L|K$. 
    According to \cite[p.\ 190]{Nk} the restriction map $G(L|K) \onto 
    G(L'|K)$ throws $G^n(L,\p)$ onto $G^n(L',\p)$. 
    Together with the previous proposition we obtain the following 
    equivalence of conditions on $L'$:
    \begin{eqnarray*}
      L'\subseteq L^n(\p) \iff G(L|L') \supseteq G^n(L,\p) 
        \iff G^n(L',\p)=1 \iff \\
      n\geq f(L',\p) \iff \m\setminus\{\p\}+n\p\geq\f(L'|K) 
        \iff L' \subseteq L \cap K_S^{\m\setminus\{\p\}+n\p}.
    \end{eqnarray*}
    The proof for the decomposition field $L^{-1}(\p)$ is similar.
  \item For a place $\p$ of $K$ with extension $\q$ to $L$ denote by 
    $d(\q|\p)$ the different exponent and by $f(\q|\p)$ the inertia degree 
    of $\q|\p$, which should not be confused with the conductor exponent 
    $f(L,\p)$. 
    Using the Hasse-Arf Theorem one can rewrite Hilbert's (lower index) 
    different formula in upper index form (for details see the appendix):
    \[
      d(\q|\p) 
               = \sum_{n=0}^{m_\p-1} \left( 
                   [L:L^0(\p)] - [L^n(\p):L^0(\p)] 
                 \right).
    \]
    The discriminant is the norm of the different. 
    Hence $\d(L|K) = \sum_\p d(L,\p)\p$ with
    \[
      d(L,\p) = \sum_{\q|\p} f(\q|\p) d(\q|\p)  
              = [L:K] m_\p - \sum_{n=0}^{m_\p-1} [L^n(\p):K],
    \]
    again using that $L^0(\p)$ is the inertia field of $\p$ in $L$.
  \qed
  \end{subproof}
  \noqed
\end{proof}

Thus computation of the discriminant (and thereby of the genus) of ray class 
field extensions amounts to determining their degrees. 
This has already been seen by Cohen et al.\ \cite{Co} in the number 
field case. 
As an illustration we give the following 

\begin{hayesex}\label{hayesex}
  Assume that $S$ consists of exactly one place of degree $d\in\N$. 
  Let $\m = \sum_\p m_\p\p > \o$ be an $S$-cycle, $h$ the (divisor) 
  class number and $g$ the genus of $K$. 
  Putting $\phi(\m) := (\I_S^\o:\I_S^\m) = \prod_{m_\p>0} (q^{\deg\p}-1) 
  q^{(m_\p-1)\deg\p}$ we have: 
  \begin{substate}
  \item
    $[K_S^\o:K]=h\cdot d$ and $[K_S^\m:K]=h\cdot d\cdot\phi(\m)/(q-1)$.
  \item The genus of $K_S^\o$ is $g(K_S^\o)=1+h(g-1)$. 
    That of $K_S^\m$ is 
    \[
      g(K_S^\m) = 1 + h\cdot[\phi(\m)(2g-2+\deg\m)-s]/(2q-2)
    \]
    with $s = (\phi(\m)/\phi(\p)+q-2)\deg\p$ if $\m$ is the multiple of a 
    single place $\p$ and $s = \sum_{m_\p>0} \phi(\m)\deg\p/\phi(\p)$ 
    otherwise. 
  \end{substate}
\end{hayesex}
  
\begin{proof}
  In order to verify the formula for $\phi(\m)$ look at the proof of 
  Proposition~\ref{idelprop} and recall that $(U_\p:U_\p^{(n)}) = 
  (q^{\deg\p}-1)q^{(n-1)\deg\p}$ for each place $\p$ and $n\in\N$.
  \begin{subproof}
    \newcommand{\s}{\mathfrak s}
    \item $|S|=1$ implies $h_S=h\cdot d$ and $\O_S^*=\F_q^*$, which is 
      mapped injectively into $\I_S^\o/\I_S^\m$ because $\m>\o$.
      Hence the degree formulas follow from the exact sequence in 
      Proposition~\ref{idelprop} 
    \item If $\deg\m=1$ then $\f(K_S^\m|K)=\d(K_S^\m|K)=\o$ by (a).
      Let $\deg\m>1$ and $m_\p>0$. 
      Then 
      \begin{eqnarray*}
        \sum_{n=1}^{m_\p-1} [K_S^{\m\setminus\{\p\}+n\p}:K]
        &=& \frac{h\cdot d}{q-1}
            \sum_{n=1}^{m_\p-1}\phi(\m) q^{(n-m_\p)\deg\p} \\
        &=& \frac{h\cdot d}{q-1}
            \frac{\phi(\m)}{\phi(\p)}(1-q^{(1-m_\p)\deg\p}),
      \end{eqnarray*}
      and $[K_S^{\m\setminus\{\p\}}:K]$ is $h\cdot d$ or 
      $\frac{h\cdot d}{q-1} \frac{\phi(\m)}{\phi(m_\p\p)}$ depending on 
      whether or not $\m$ is a multiple of $\p$. 
      Applying~\ref{raythm}, all this results in 
      $\d(K_S^\m|K) = hd(\phi(\m)\m-\s)/(q-1)$ with 
      $\s = (\phi(\m)/\phi(\p)+q-2)\p$ if $\m$ is supported by merely one 
      place $\p$ and $\s = \sum_{m_\p>0} \frac{\phi(\m)}{\phi(\p)} \p$ 
      otherwise. 
      The genus is now easily computed from the degree and the discriminant 
      by means of the Hurwitz genus formula. 
      (See \cite[p.\ 88]{St} and take into consideration that 
      $K_S^\m$ has constant field $\F_{q^d}$.)
    \qed
  \end{subproof}
  \noqed
\end{proof}

As soon as $S$ consists of more than one place, the determination of the 
degrees is much more difficult. 
We shall dedicate Section 3 to working out this problem in a special case. 

%% file: bounds.tex
\section{Upper bounds}

As before let $K$ be a global function field with full constant field 
$\F_q$. 
By the Hasse-Weil Theorem (see \cite[p.\ 169f]{St}) as a first estimate 
for its number $N(K)$ of rational places we obtain 
\begin{equation}\label{weilbound}
  N(K) \leq q + 1 + 2g(K) \sqrt q,
\end{equation}
where $g(K)$ denotes the genus of $K$. 
This is called the Hasse-Weil bound. 
In case equality holds, $K$ is called \df{maximal}. 
Obviously $K$ can only be maximal if $g(K)=0$ or if $q$ is a square. 
A consideration involving the constant field extension of degree 2 (see 
\cite[p.\ 182]{St}) and a recent work by Fuhrmann and Torres \cite{FT} 
shows that in addition one must have $g(K) = (q-\sqrt q)/2$ or 
$g(K) \leq (\sqrt q - 1)^2/4$. 
In the non-square case Serre \cite{S1} was able to 
improve~\eqref{weilbound} by 
\[
  N(K) \leq q + 1 + g(K) \lfloor 2\sqrt q \rfloor
\]
where $\lfloor a \rfloor$ denotes the integer part of $a\in\R$.

The maximum number of rational places a global function field of genus $g$ 
with full constant field $\F_q$ can have is usually denoted $N_q(g)$.
We call $K$ \df{optimal} if $N(K) = N_q(g(K))$.
Serre \cite{S3} succeeded in determining $N_q(g)$ for $g=1$ and 2 in 
general and in some cases also for $g=3$ and 4.
An overview of what is known about $N_q(g)$ for $g\leq50$ and $q$ a not 
too large power of 2 or 3 is found in \cite{GV7}, \cite{Ni}, \cite{NX11} or \cite{NX12}.

For large genera one can improve~\eqref{weilbound} by means of the so 
called explicit Weil formulas (see \cite[p.\ 183]{St}):
Any finite sequence $(c_n)$ of reals $c_n\geq0$ satisfying
\[
  \textstyle
  1 + \sum_{n\geq1} 2 c_n \cos n\ph \geq 0 \quad \forall \ph\in \mathbb R
\]
provides an estimate 
\begin{equation}\label{explifor}
  \textstyle
  g(K) \geq (N(K)-1)\sum_{n\geq1} c_n q^{-n/2} - \sum_{n\geq1} c_n q^{n/2}.
\end{equation}
We shall give a brief outline of the maximization of the right hand side due 
to  Oesterl\'e as found in \cite{Sf} and, based on it, derive an estimate 
for $S$-class numbers that spares us from going into the explicit formulas 
again and again. 

Given a real number $N \geq q+1$, let $m$ be the integer $\geq2$ such that 
$q^{m/2} \leq N-1 < q^{(m+1)/2}$ and $u := 
[(q^{(m+1)/2}-q^{(m-1)/2})/(N-1-q^{(m-1)/2})-1] / \sqrt q \in (0,1]$.
Because the continuous maps
\begin{map}{F_m}
  [\frac\pi{m+1},\frac\pi m] & \to & [0,1] \\
  \ph & \mapsto & - \left({\cos\frac{m+1}2\ph}\right)
        \bigm/ \left({\cos\frac{m-1}2\ph}\right)
\end{map}
are onto and strictly increasing, $\th_q(N) := \cos(F_m^{-1}(u))$ defines 
a continuous bijection $\th_q : [q+1,\infty) \to [0,1)$, which is also 
increasing. 
\begin{oesthm}\label{oesthm}
  Define
  \begin{map}{g_q}
    [q+1,\infty) & \to & [0,\infty) \\
    N & \mapsto & 1 + \frac{(\sqrt{q}\th_q(N)-1)N}{q-2\sqrt{q}\th_q(N)+1}. 
  \end{map}
  and assume that $N(K) \geq q+1$. 
  Then $g(K) \geq g_q(N(K))$, and $g_q(N(K))$ is the maximum of the right 
  hand side of~\eqref{explifor} for $q\geq3$ (but not always for $q=2$).
\end{oesthm}
After checking that $g_q$ is strictly increasing we find that the value at 
$g\in\N_0$ of its inverse $\bar N_q := g_q^{-1}$ provides an upper bound 
for the number $N_q(g)$, which for $g \leq (q - \sqrt q)/2$ coincides with 
the Hasse-Weil bound~\eqref{weilbound} and improves it for larger $g$. 
$\bar N_q(g)$ will be called the Oesterl\'e bound.

Now let $S$ be a non-empty set of places of $K$.
Applying Oesterl\'e's Theorem to the Hilbert class field of $\O_S$ yields 
an upper bound for the $S$-class number $h_S$.
\begin{lem}\label{hbarlem}
  Assume that $K$ has genus $g\geq1$ and that $S$ contains at least 
  $N>(\sqrt q-1)(g-1)$ rational places, then 
  \[
    h_S \leq \hbar_q(g,N) := \th_q^{-1}(t)/N \quad\mbox{with}\quad
      t := (1 + \frac{(q-1)(g-1)}{N+2g-2}) / \sqrt q.
  \]
\end{lem}
\begin{proof}
  The restrictions on $g$ and $N$ ensure that $1/\sqrt q \leq t < 1$, whereby 
  $\hbar_q$ is well-defined. 
  Let $g' := g(K_S^{\o})$ and $N' := N(K_S^{\o})$ be the genus \resp the 
  number of rational places of the Hilbert class field of $\O_S$, which
  satisfy $h_S(g-1) = g'-1$ and $h_SN \leq N'$. 
  For $N' \leq q+1$ the claimed estimate is obvious. 
  Assume $N' \geq q+1$, then by Oesterl\'e's Theorem 
  \[
    \frac{g-1}N \geq \frac{g'-1}{N'} \geq \frac{g_q(N')-1}{N'}
    = \frac{\sqrt q\th_q(N')-1}{q-2\sqrt q\th_q(N')+1}.
  \]
  Isolating $\th_q(N')$ in this inequality results in $\th_q(N') \leq t$, 
  from which by the isotony of $\th_q$ we conclude the assertion. 
\end{proof}
Note that $\th_q^{-1}$ can be computed very efficiently using the cosine 
formula for multiple angles. 
In fact, on each interval $[\cos\frac\pi m,\cos\frac\pi{m+1}]$ with $m$ 
being an integer $\geq2$ it is identical with a $\Q(\sqrt q)$-rational 
function, namely $\th_q^{-1}(t) = 1 + (q^{(m+1)/2}+q^{m/2}u)/(1+\sqrt qu)$ with $u = f_{m+1}(t)/f_{m-1}(t)$, where
\[
  f_n(t) := \sum_{i=0}^{\floor{n/2}} \tfrac n{n-i} \tbinom{n-i}i 
            (-2t-2)^{\floor{n/2}-i} \text{ for } n\in\N.
\]
This also shows $\hbar_q(g,N) \in \Q(\sqrt q)$.


An essential step in determining the degrees of the $S$-ray class fields 
is the computation of $h_S$. This can often be done using Lemma~\ref{hbarlem} and the relation of $h_S$ with the \df{$S$-regulator} $\reg_S$, defined as the index of the group $(\O_S^*)$ of principal divisors of $S$-units in the group $\D_S^0$ of divisors of degree 0 with support in $S$ (see \cite{Ro1} for details). 

\begin{ex}\label{ex020205}
  Consider the function field $K := \F_2(x,y)$ with $y^2+y = x^3(x+1)^2$ 
  of genus $g(K)=2$.
  It has 5 places of degree 1 and no place of degree 2, whence its 
  divisor class number is $h(K)=13$. 
  Let $\p_\infty$ be the pole of $x$ (and $y$) and put 
  $\p_0:=(0,0)$, $\p_1:=(0,1)$, $\p_2:=(1,0)$ and $\p_3:=(1,1)$, where 
  $(\alpha,\beta)$ denotes the common zero of $x-\alpha$ and $y-\beta$.
  First look at the set $S := \{ \p_\infty, \thru\p03 \}$ of all rational 
  places of $K$. 
  Since $\hbar_2(2,5)=(2503+960\sqrt2)/2911 < 2$ and 
  $h_S \cdot \reg_S = h(K)$ we see that $h_S=1$ and $\reg_S=13$. 
  Let $O$ be the subgroup of $\O_S^*$ generated by $x$, $x+1$, $y$ and 
  $y+x^2$. 
  Without regard to sign, the index of $(O)$ in $\D_S^0$ equals any   
  $4\times4$-minor of the matrix 
  \[
    D := \left( 
      v_{\p_j}(z) 
    \right)_{\on{z=x,x+1,y,y+x^2}{j=\infty,1,2,3,0}} = 
    \begin{array}({ccccc})
      -2&1&0&0&1\\
      -2&0&1&1&0\\
      -5&0&2&0&3\\
      -5&0&0&3&2
    \end{array},
  \]
  which turns out to be $\pm13$. 
  Hence $\O_S^*=O$. 
  Now we can compute $h_{S'}$ and a basis of $\O_{S'}^*$ for any 
  $S' \subseteq S$, \eg for $S_r := \{ \p_\infty,\thru\p1r \}$ the first 
  $r \in \{1,2,3\}$ rows of the Hermite normal form 
  \[
    \begin{array}({rrrr})
      13&6&-3&-2\\
       5&3&-1&-1\\
      -5&-2&1&1\\
       1&0&0&0
    \end{array} \cdot D = 
    \begin{array}({rrrrr})
      -13&13&0&0&0\\
       -6& 5&1&0&0\\
        4&-5&0&1&0\\
       -2& 1&0&0&1
    \end{array}
  \]
  considered as elements of $\D_S^0$ constitute a basis of $(\O_{S_r}^*)$, 
  so we see that $h_{S_r}=1$.
  This example is continued in Part (B) of Section 3.
\end{ex}

%% file: many.tex
\section{Many rational places}

\newcommand\tablesize\footnotesize

\newenvironment {flattable}[2]
{
\tablesize\begin{tabular}{*{#2}{|c}|}
  \hline
  \multicolumn{#2}{|c|}{$q=#1$} \\
  \hline
}
{ \hline \end{tabular}}

\newenvironment {flattable*}[1]
{
\tablesize\begin{tabular}{*{#1}{|c}|}
  \hline
}
{ \hline \end{tabular}}

\newenvironment {bigtable} [1]
{\tablesize\begin{tabular}[t]{|r|r|c|c|c|c|c|}
  \hline
  \multicolumn{7}{|c|}{$q=#1$} \\
  \hline
  $g$ & $N_q(g)$ & $n$ & $l$ & $|S|$ & $h_S$ & $g(K)$ \\
  \hline
}
{ \hline\end{tabular} }

\newenvironment {slimtable} [1]
{\tablesize\begin{tabular}[t]{|r|r|c|c|c|}
  \hline
  \multicolumn{5}{|c|}{$q=#1$} \\
  \hline
  $g$ & $N_q(g)$ & $n$ & $l$ & $|S|$ \\
  \hline
}
{ \hline\end{tabular} }

\newcommand{\mergetable}[1]
{
  \hline
  \multicolumn{5}{c}{} \\
  \hline
  \multicolumn{5}{|c|}{$q=#1$} \\
  \hline
  $g$ & $N_q(g)$ & $n$ & $l$ & $|S|$ \\
  \hline
}

When searching for curves with many rational points, it seems to be 
fruitful to investigate ray class fields $K_S^\m$, where $K$ already has 
many rational places compared to its genus, several or all of which 
are contained in $S$. 

First we want to fix some notation. 
Let $K|\F_q$ be as before, and let $g(K)$ be the genus and $p$ the 
characteristic of $K$. 
Then $q=p^e$ for some $e\in\N$. 
By $v_p$ we denote the (normalized) $p$-adic valuation on $\Q$ 
or its $p$-adic completion $\Q_p$ with values in $\Z\cup\{\infty\}$. 
Throughout this whole section, $n$ always denotes a positive (or, where this makes sense, a non-negative) integer.
For $n\geq1$ we put $n^* := n/p^{v_p(n)}$. 

To simplify the situation somewhat we restrict ourselves to the case 
of $\m$ being the multiple of a single rational place $\p$ of $K$.\footnote{
  The more general case of $\p$ having arbitrary degree is treated 
  in~\cite{Au}
} 
Let $S$ be a non-empty set of places of $K$ not containing $\p$. 
Put $U^{(n)} := U_\p^{(n)}$ and $U_S := \O_S^* \cap U^{(1)}$. 
By Dirichlet's Unit Theorem $U_S \simeq \O_S^*/\F_q^*$ is a free abelian 
group of rank $|S|-1$. 
Since $U^{(0)}=\F_q^*U^{(1)}$ and $\F_q^*$ is contained in $\O_S^*$, the 
exact sequence of Proposition~\ref{idelprop} yields $K_S^\o=K_S^\p$, and using the 
isomorphisms discussed in Section~1, it can be rewritten as 
\begin{equation}\label{useq}
  1 \to U_S\cap U^{(n)} \to U_S \to U^{(1)}/U^{(n)} \to G(K_S^{n\p}|K) 
    \to \Cl(\O_S) \to 1,
\end{equation}
whence $G(K_S^{n\p}|K_S^\o) \simeq U^{(1)}/U_SU^{(n)}$.
Let $\la_S^{(n)} := \log_p [K_S^{n\p}:K_S^\o]$, where $\log_p$ means 
logarithm to the base $p$. 
Then clearly $\la_S^{(0)} = \la_S^{(1)} = 0$ and 
$\la_S^{(n)} = \log_p(U^{(1)}:U_S U^{(n)}) \in \N_0$ for $n\geq1$.
Wanting to have many examples, for each $l \in 
\set{\la_S^{(n-1)}}{\la_S^{(n)}}$ by Galois theory we choose an 
intermediate field $L_{l,S}$ of $K_S^{n\p}|K_S^{(n-1)\p}$ of degree $p^l$ 
over $K_S^{\o}$.
Assume $\deg S = 1$, then by Proposition~\ref{rayprop} each $L_{l,S}$ has again full 
constant field $\F_q$ and Theorem~\ref{raythm} together with the Hurwitz genus 
formula yields 
\begin{equation}\label{genfor}
  g(L_{l,S}) = 1 + \frac{h_S}2 \left( p^l(2g(K)-2+n) - 
             \sum_{\nu=0}^{n-1} p^{\la_S^{(\nu)}}
           \right).
\end{equation}
Since $h_S/h_{S\cup\{\p\}}$ is the inertia degree of $\p$ in $L_{l,S}$, 
the number $N(L_{l,S})$ of rational places of $L_{l,S}$ satisfies 
\begin{equation}\label{numfor}
  N(L_{l,S}) \geq h_S \cdot p^l \cdot |\{ \q\in S \sothat \deg\q=1 \}| 
  + \left\{\begin{array}{cl}
      h_S & \mbox{if } h_S=h_{S\cup\{\p\}}, \\
      0   & \mbox{if } h_S>h_{S\cup\{\p\}},
    \end{array}\right.
\end{equation}
where equality holds if and only if $L_{l,S} \not\subseteq 
K_{S\cup\{\q\}}^{n\p}$ 
(which is guaranteed for $l>\la_{S\cup\{\q\}}^{(n)}$) 
for each rational place $\q$ not contained in $S\cup\{\p\}$. 

We shall give two methods to determine the degrees of the $K_S^{n\p}$ over $K$.
Method (A) is specialized to the rational function field and can be carried out very easily in practice. 
Method (B) means much more effort but applies to arbitrary $K$.

\input{rational.tex}
\input{method.tex}
\input{tables.tex}

%% file: rational.tex
\subsection{Rational function field}

First we investigate the case $K=\F_q(x)$, where $x$ is an indeterminate 
over $\F_q$. 
Also let $S$ consist of rational places only (implying $h_S=h(K)=1$, so we have $K_S^\o = K$). 
For $\alpha\in\F_q$ denote the zero of $x-\alpha$ by $\p_\alpha$. 
After an appropriate transformation of the variable $x$ we can assume that $\p$ is the pole of $x$ and that $\p_0$ is contained in $S$. 
Then $\pi := 1/x$ is a uniformizer at $\p$, and $U_S$ is freely generated by the functions $1 - \alpha \pi$ with 
$\alpha \in A_S := \{ \alpha\in\F_q^* \sothat \p_\alpha\in S \}$. 

By what has been said above, our task is to determine the logarithmic 
degrees $\la_S^{(n)} = \log_p[K_S^{n\p}:K]$.
In many cases they obey a very simple rule (Theorem~\ref{ratthm} below), which can be developed by the following idea involving Newton's Formulas. 

Choose a generator $\om$ of $\F_q^*$, put $I := \set1{q-1}$ and
$I_S := \{ j \in I \sothat \om^j \in A_S \}$, let 
$\bar{\makebox[.5em]{}}$ denote the canonical projection 
\[
  \F_p[x] \onto R := \F_p[x]/(x^{q-1}-1)
\]
and define the morphism 
\begin{map}{\rho}
  U_S & \to & R \\
  \prod_{\alpha\in A_S} (1 - \alpha\pi)^{c_\alpha} & \mapsto &
  \sum_{j\in I_S} c_{\om^j} \bar x^j
\end{map}
of groups with kernel $U_S^p$ and image 
$R_S := \bigoplus_{j\in I_S} \F_p \bar x^j $. 

Denote by $f_\alpha \in \F_p[x]$ the minimal polynomial of 
$\alpha\in\F_q$ over $\F_p$. 
The ring $R$ is frequently used in coding theory. 
It is principal (as homomorphic image of $\F_p[x]$) and its ideals (which are in fact cyclic codes) correspond bijectively to the monic divisors of $x^{q-1}-1$. 
Hence the kernel of the $\F_p$-algebra morphism 
\begin{map}{\Phi^{(n)}}
  R & \to & \F_q^{n-1} \\
  \bar f & \mapsto & (f(\om^i))_{1\leq i<n}
\end{map}
is just the ideal $R \bar f^{(n)}$, where $f^{(n)}$ is the least common 
multiple of all $f_{\om^i}$ with $1\leq i<n$. 

Denote the degree of $f^{(n+1)}/f^{(n)}$ by $e^{(n)}$ and call $n$ 
\df{(an) initial (of $q$)} if $e^{(n)}>0$. 
Let $G(\F_q|\F_p)$ act on $I$ via the bijection $I\leftrightarrow\F_q^*$, 
$i\mapsto\om^i$, and consider the orbits under this action.
Obviously $n$ is initial iff it is in $I$ and is the least element in its orbit, in which case $e^{(n)} = \deg f_{\om^n} = [\F_p(\om^n):\F_p]$ is the length of this orbit and divides $e$. 
For example if $q=16$, $I$ decomposes into the orbits $\{1,2,4,8\}$, $\{3,6,9,12\}$, $\{5,10\}$, $\{7,11,13,14\}$ and $\{15\}$, which means that the initials of $q=16$ are 1, 3, 5, 7 and 15, where 
$e^{(1)} = e^{(3)} = e^{(7)} = 4$, $e^{(5)} = 2$ and $e^{(15)} = 1$.

Let $\ph$ be the Frobenius on $K$, which takes everything to its $p$-th power. 
Writing the elements of $I$ in $p$-adic representation with $e$ digits, the action of $\ph^l$ is just a rotation of these digits by $l$ places. 
Using this to determine whether or not $n$ is initial, the proof of the following lemma is easy and can therefore be omitted. 
\begin{lem}\label{eflem}
  \begin{substate}
  \item $f^{(n)} = x^{q-1}-1$ iff $n \geq q$.
  \item $f^{(n)} = 1+\cdots+x^{q-2}$ iff $q-q/p \leq n < q$.
  \item If $p$ divides $n$ then $e^{(n)}=0$.
  \item If $e$ is even and $n \in \set1{2\sqrt q} \setminus 
    (\Z p \cup\{\sqrt q+1\})$, then $e^{(n)}=e$, $e^{(\sqrt q+1)}=e/2$
    and $e^{(2\sqrt q+1)}=0$.
  \item If $e$ is odd, then $e^{(n)}=e$ for $n \in \set1{\sqrt{pq}} 
    \setminus \Z p $ and $e^{(n)}=0$ for $n \in$ 
    \mbox{$\set{\sqrt{pq}}{\sqrt{pq}+p}$}. 
  \end{substate}
\end{lem}

Put $n' := \lceil n/p \rceil$, where $\lceil u \rceil$ denotes the least 
integer greater or equal to $u$. 
Then the Frobenius map induces an injection $\ph^{(n)} : U^{(1)}/U^{(n')} 
\into  U^{(1)}/U^{(n)}$. 
To be able to continue we must impose certain restrictions on $n'$ and 
$n$ in terms of the two numbers 
\begin{eqnarray*}
  n_S' & := & \max \{ n\in\N \sothat \la_S^{(n)}=0 \} = 
             \max \{ n\in\N \sothat U^{(1)}=U_SU^{(n)} \} 
             \quad\mbox{and} \\ 
  n_S & := & \min \{ q-q/p, n\in\N \sothat R_S \cap R \bar f^{(n)} = 0 \} 
             - 1.
\end{eqnarray*}
The main work will be done in the following
\begin{lem}\label{ratlem}
  \begin{substate}
  \item Let $y\in U_S$. Then $y U^{(n)}$ lies in the image of $\ph^{(n)}$ 
    if and only if $\rho(y) \in R \bar f^{(n)}$. 
  \item Assume that $n' \leq n_S'$ or $n > n_S$. Then
    \[
      \frac{(U_S U^{(n)}:U^{(n)})}{(U_S U^{(n')}:U^{(n')})} 
      = (R_S:R_S \cap R \bar f^{(n)}).
    \]
  \end{substate}
\end{lem}

\begin{proof}
  \begin{subproof}
  \item Write $y=\prod_{\alpha\in A_S} (1-\alpha\pi)^{c_\alpha}$ with 
    integers $c_\alpha$, let
    \[
      y = \sum_{j=0}^\infty \sigma_j \pi^j, \quad \sigma_j\in\F_q,
    \]
    be the expansion of $y$ at $\p$ and put 
    $s_i := \sum_{\alpha \in A_S} c_\alpha \alpha^i$.
    Since multiplication of $y$ by $(1-\alpha\pi)^{p^l}$ neither changes 
    the $s_i$ nor $\thru\sigma0n$ provided $n<p^l$ and $l\in\N$, 
    we can assume that none of the $c_\alpha$ is negative and obtain
    \[
      n\sigma_n + \sum_{i=1}^n s_i \sigma_{n-i} = 0
    \]
    from Newton's Formulas.
    By induction we deduce the equivalence
    \[
      \forall j\in\set1{n-1}\setminus \Z p  : \sigma_j=0
      \iff s_1=\ldots=s_{n-1}=0.
    \]
    Since $\Phi^{(n)}(\rho(y)) = (s_i)_{1\leq i<n}$, this can be 
    reformulated to yield the assertion.
  \item All the above defined morphisms can be arranged in a commutative 
    diagram
    \[\begin{array}{rcccccccl}
        &     &  1         &&  1         && 0 \\
        &     & \downarrow && \downarrow && \downarrow \\
      1 & \to & U_S \cap U^{(n')} & \stackrel\ph\to & U_S \cap U^{(n)}
        & \stackrel{\rho}\to & R_S \cap R \bar f^{(n)} & \to & 0 \\
        &     & \downarrow && \downarrow && \downarrow \\
      1 & \to & U_S & \stackrel\ph\to & U_S
        & \stackrel{\rho}\to & R_S & \to & 0 \\
        &     & \downarrow && \downarrow 
        && \makebox[4.1em][r]{$\downarrow{\scriptstyle\Phi^{(n)}}$} \\
      1 & \to & U^{(1)}/U^{(n')} & \stackrel{\ph^{(n)}}\to 
        & U^{(1)}/U^{(n)} && \F_q^{n-1} \\
    \end{array}\]
    (use (a) to see that $\rho(U_S \cap U^{(n)}) \subseteq R \bar 
    f^{(n)}$), the exactness of which is easily verified except at one 
    critical point, namely the horizontal passage through the upper right 
    corner, \ie it remains to show that $\rho(U_S \cap U^{(n)})$ really 
    is all of $R_S \cap R \bar f^{(n)}$. 
    To this end let $y \in U_S$ with $\rho(y) \in R \bar f^{(n)}$.

    First assume that $n' \leq n_S'$. Then by (a) there exists $y'\in 
    U_S$ with $yU^{(n)} = \ph^{(n)}(y'U^{(n')}) = {y'}^p U^{(n)}$, \ie 
    $y/{y'}^p \in U^{(n)}$. But $\rho(y/{y'}^p) = \rho(y)$ from the exactness 
    of the middle row. 
    
    By Lemma~\ref{eflem}(a)--(b), $R_S \cap R \bar f^{(n)}$ is trivial if 
    $n \geq q$ and also if $n>n_S$ and $|S|<q$. 
    Finally, if $q-q/p\leq n<q$ and $|S|=q$, from~\ref{eflem}(b) we know that 
    $\rho(y) = \gamma(1+\cdots+x^{q-2})$ for some $\gamma\in\F_p$, \ie 
    $y$ is a power of $1-\pi^{q-1}$ and therefore contained in $U^{(n)}$.

    Now use the isomorphism
    \[
      U_SU^{(n)}/U^{(n)} \simeq \frac {U_S/\ph(U_S \cap U^{(n')})} 
                         {U_S \cap U^{(n)}/\ph(U_S \cap U^{(n')})}
    \]
    and the exactness just proven to obtain
    \[
      \frac {(U_SU^{(n)}:U^{(n)})} {(U_SU^{(n')}:U^{(n')})}
      = \frac { (U_S:\ph(U_S)) (\ph(U_S):\ph(U_S \cap U^{(n')})) }
              { (U_S \cap U^{(n)}:\ph(U_S \cap U^{(n')}))
                (U_S:U_S \cap U^{(n')}) } =
    \]
    $= (R_S:R_S \cap R \bar f^{(n)})$.
  \qed
  \end{subproof}
  \noqed
\end{proof}

We have arrived at the main theorem of Method~(A). 
It gives a formula for the (logarithmic) degrees of $K_S^{n\p}|K$ in 
terms of the non-$p$-part $n^*$ of $n$ defined above and the numbers
\[
  e_S^{(n)} := \dim_{\F_p} \left( R_S \cap R \bar f^{(n)} 
                         / R_S \cap R \bar f^{(n+1)} \right),
\]
which are $\leq e^{(n)}$ and whose computation is an easy exercise in 
linear algebra or even simpler. 
Lauter \cite{Lt3} has already proved this formula for all $n$ in 
case $|S|=q$ by means of a different method.
\begin{thm}\label{ratthm}
  \begin{substate}
  \item As long as $n<pn_S'$ the $\la_S^{(n)}$ satisfy the recursion
    \[
      \la_S^{(n+1)} = \la_S^{(n)} + e - e_S^{(n^*)}.
    \]
    If $pn_S' \geq n_S$ this formula holds for every $n\in\N$.
  \item $n_S' = \min \{ n\in\N \sothat e_S^{(n^*)}<e \}$ and
    \[
      n_S = \left\{\begin{array}{cl}
              0 & \mbox{if } |S|=1, \\
              \max \{ n\in\N \sothat e_S^{(n)}>0 \} 
                 & \mbox{if } 1<|S|<q, \\
              q-q/p-1 & \mbox{if } |S|=q.
            \end{array}\right.
    \]
  \item If $e=1$, \ie $q=p$, then 
    \[
      e_S^{(n)} = \left\{\begin{array}{cl}
        1 & \mbox{if } n<|S|, \\
        0 & \mbox{otherwise}
      \end{array}\right.
    \]
    and the formula in (a) holds for every $n\in\N$.
  \end{substate}
\end{thm}

\begin{proof}
  \begin{subproof}
  \item From the exact sequence~\eqref{useq} we obtain 
    $\la_S^{(n)} = (n-1)e - \log_p(U_SU^{(n)}:U^{(n)})$. 
    Hence it suffices to show that
    \[
      \log_p \frac {(U_SU^{(n+1)}:U^{(n+1)})} {(U_SU^{(n)}:U^{(n)})}
      = e_S^{(n^*)}
    \]
    under the given assumptions.
    But this is easily derived from the previous lemma by induction on 
    $n$ while distinguishing whether or not $p$ divides $n$.
  \item follows immediately from (a).
  \item Observe that
    \[
      \dim_{\F_p} \Phi^{(n)}(R_S) 
      = \rank (\om^{ij})_{\on{1\leq i<n}{j\in I_S}}
      = \min\{n-1,|I_S|\},
    \]
    showing the assertion on $e_S^{(n)}$,
    and use (b).
  \qed
  \end{subproof}
  \noqed
\end{proof}

Without any further effort, Theorem~\ref{ratthm} (a) and (c) in combination with equations~\eqref{genfor} and~\eqref{numfor} leads to the following tables, giving $g(L_{l,S})$ and $N(L_{l,S})$ for some prime numbers $q=p$ and some values of $l$ and $|S|$. 
\begin{figure}
\begin{flattable}2{12}
  $l,|S|$      &1,2&2,2&3,2 &4,1 &4,2 &5,1 &5,2 &6,1 &6,2 &7,2 & 8,2 \\ 
  \hline
  $g(L_{l,S})$ & 1 & 5 & 15 & 17 & 39 & 49 &103 &129 &247 &567 &1271 \\
  \hline
  $N(L_{l,S})$ & 5 & 9 & 17 & 17 & 33 & 33 & 65 & 65 &129 &257 & 513 \\
\end{flattable}
\\[.5\tablesep]
\begin{flattable}3{12}
  $l,|S|$      &1,2&1,3&2,2 &2,3 &3,1 &3,2 &3,3 &4,2 &4,3 &5,2 &5,3  \\ 
  \hline
  $g(L_{l,S})$ & 1 & 3 & 10 & 15 & 21 & 46 & 69 &181 &258 &667 & 987 \\
  \hline
  $N(L_{l,S})$ & 7 &10 & 19 & 28 & 28 & 55 & 82 &163 &244 &487 & 730 \\
\end{flattable}
\\[.5\tablesep]
\begin{flattable}5{12}
  $l,|S|$      &1,2&1,3&1,4&1,5&2,2&2,3&2,4 &2,5 &3,2 &3,3 &3,4 \\ 
  \hline
  $g(L_{l,S})$ & 2 & 4 & 6 &10 &22 &34 & 56 & 70 &172 &284 &356 \\
  \hline
  $N(L_{l,S})$ &11 &16 &21 &26 &51 &76 &101 &126 &251 &376 &501 \\
\end{flattable}
\\[.5\tablesep]
\begin{flattable}7{12}
  $l,|S|$      &1,2&1,3&1,4&1,5&1,6&2,2&2,3 &2,4 &2,5 &2,6 &2,7 \\ 
  \hline
  $g(L_{l,S})$ & 3 & 6 & 9 &12 &15 &45 & 69 & 93 &117 &162 &189 \\
  \hline
  $N(L_{l,S})$ &15 &22 &29 &36 &43 &99 &148 &197 &246 &295 &344 \\
\end{flattable}
\caption{Ray class fields of $\F_q(x)$ for $q=2$, 3, 5 and 7}
\end{figure}
For $q=2$ examples from \cite{XN} are reproduced, and all entries with 
$|S|=1$ are known as cyclotomic function fields. 
For $q=3$ and $g(L_{l,S})\geq181$ we encounter first examples of new function 
fields with many rational places. 
The table for $q=5$ continues a series of examples worked out 
in~\cite{NX8} for genus up to $56$. 
For $q=7$ all examples except the one with $g(L_{l,S})=3$ are new. 
Writing down the tables for larger prime numbers is not worthwhile 
because the entries are getting too bad. 

A direct computation shows that $1$ and $3$ are the initials of $q=4$
with $e_S^{(1)}=\min\{2,|S|-1\}$ and $e_S^{(3)}=\max\{0,|S|-3\}$, so
Theorem~\ref{ratthm} gives the table below. 
It yields some new estimates compared to~\cite{NX12}, \eg $97 \leq N_4(67) \leq 117$. 
\begin{figure}
\begin{flattable}4{12}
  $l,|S|$      &1,4&3,1&2,3&2,4&4,1&3,3 &3,4 &5,1 &4,3 &4,4 &5,2 \\ 
  \hline
  $g(L_{l,S})$ & 1 & 2 & 3 & 5 & 6 & 11 & 13 & 22 & 27 & 33 & 37 \\
  \hline
  $N(L_{l,S})$ & 9 & 9 &13 &17 &17 & 25 & 33 & 33 & 49 & 65 & 65 \\
\end{flattable}
\nopagebreak
\\
\begin{flattable*}{11}
  $l,|S|$      &5,3 &5,4 &6,2 &6,3 &6,4 &7,2 &7,3 &7,4 &8,3 & 8,4 \\
  \hline
  $g(L_{l,S})$ & 67 & 81 &101 &147 &177 &229 &339 &433 &723 & 945 \\
  \hline
  $N(L_{l,S})$ & 97 &129 &129 &193 &257 &257 &385 &513 &769 &1025 \\
\end{flattable*}
\caption{Ray class fields of $\F_4(x)$, and intermediate fields}
\end{figure}
A more complete table obtained by various other constructions is found \eg in~\cite{NX4}. 

While in the previous examples the $e_S^{(n)}$ only depended on $|S|$, they may vary for different $S$ in the same cardinality when regarding proper prime powers $q\neq4$. 
Due to this fact we have to perform some linear algebra before computing the tables for $q=8$, 9, 16, 25, 27, 32, 49, 64 and 81, which are found at the end of this section. 

Let us now take a closer look at the situation when $|S|=q$.
Put $r:=\sqrt q$ or $\sqrt{pq}$ depending on whether $e$ is even or odd. 
From Theorem~\ref{ratthm} and Lemma~\ref{eflem}(c)--(e) we obtain two canonical 
fields, namely $K_S^{(r+2)\p}$ and $K_S^{(2r+2)\p}$ of degree $r$ and $rq$ 
over $K$ if $e$ is even, and $p-1$ canonical fields $K_S^{(r+i+1)\p}$ with 
$1\leq i<p$ of degree $q^i$ over $K$ in case $e$ is odd. 
Their genera and those of all intermediate fields $L_{l,S}$ are easily 
computed from our formulas.
\begin{cor}\label{llcor}
  Let $|S| = q$. Then $N(L_{l,S}) = 1 + p^l q$ and
  \begin{substate}
  \item $g(L_{l,S}) = \left\{\begin{array}{cl}
          \frac r2 (p^l-1) & \mbox{for } 0 \leq l \leq e/2, \\
          \frac r2 (2p^l-r-1) & \mbox{for } e/2 \leq l \leq 3e/2
         \end{array}\right\}$
    if $e$ is even.
  \item $g(L_{l,S}) = \frac12 \left( 
          p^l(r+i-1) - r - \frac{q^i-q}{q-1} 
        \right)$
    for $(i-1)e \leq l \leq ie$ with $1\leq i<p$ if $e$ is odd.
  \end{substate}
\end{cor}
$L_{l,S}$ is maximal for $0 \leq l \leq e/2$ and $e$ even.
If $e$ is odd and $p=2$ or 3 then $K_S^{(r+p)\p} = L_{(p-1)e,S}$ is optimal 
(\cf \cite{S3}). 
By comparison of conductors it is easily verified that 
$K_S^{(r+2)\p} = K(y)$ where $y$ satisfies $y^r+y = x^{r+1}$ if $e$ is 
even, and $y^q-y = x^{r/p}(x^q-x)$ if $e$ is odd. 
For the connection with the Deligne-Lusztig curves associated to the Suzuki 
and the Ree group see \cite{Lt2} and \cite{GS}. 

For odd $e$, generalizing \cite{Pd} one can write down defining equations for the first $p-1$ ray class fields.
\begin{prop}
  Let $e$ be odd and let $x_i$ satisfy $x_i^q-x_i = x^{ir/p}(x^q-x)$. 
  Then $K_S^{(r+i+1)\p} = K(\thru x1i)$ for $1\leq i<p$. 
\end{prop}
The proof is carried out in~\cite{Au} using composites of Artin-Schreier extensions. 
There you also find defining equations for other ray class fields of the rational function field. 

%% file: method.tex
\subsection{General method} 

Now we develop a method to determine the numbers $\la_S^{(n)}$ for arbitrary ground field $K$ using ideas from Xing and Niederreiter \cite{XN} and thereby generalize what they call their `first construction'. 

We fix a basis $B$ of $\F_q$ over $\F_p$ and a uniformizer 
$\pi \in K_\p^*$ at $\p$. 
For each $j \in \N$ and $\beta \in B$ consider the cyclic subgroup 
\[
  G_{j\beta}^{(n)} := \left< (1+\beta\pi^j)U^{(n)} \right>
\]
of $U^{(1)}/U^{(n)}$.
Its logarithmic order is 
\[
  \log_p \left|G_{j\beta}^{(n)}\right| = \round nj := 
  \left\{\begin{array}{cl}
     \left\lceil\log_p\frac nj\right\rceil & \mbox{if }j\leq n \\
                   0                       & \mbox{if }j\geq n.
  \end{array}\right.
\]
Denote by $\N^* := \N\setminus \Z p$ the set of all positive integers not divisble by $p$. 
The numbers $\round nj$ satisfy the easily verified properties
\begin{eqnarray}
 \label{rounditerate}
 \round{n+1}j & = & \round nj + \left\{\begin{array}{cl}
   1 & \mbox{if there is } l\in\N_0 \mbox{ such that } n = jp^l, \\
   0 & \mbox{otherwise,} \label{roundsum}
 \end{array}\right. \\
 \sum_{j\in \N^*} \round nj & = & n-1,
\end{eqnarray}
\begin{equation} \label{roundineq}
  i,j \in \N, l\in\Z, ip^l \leq j \leq n \implies  
  \round ni \geq \round nj + l.
\end{equation}
A straightforward proof using equation~\eqref{roundsum} yields the decomposition
\[
  U^{(1)}/U^{(n)} = \prod_{\on{j\in \N^*}{\beta\in B}} 
                      G^{(n)}_{j\beta}
\]
as a direct product.
Thus we obtain a non-canonical isomorphism of finite $p$-groups
\begin{map}{\mu^{(n)}}
  U^{(1)}/U^{(n)} & \to & \M^{(n)} := 
    \displaystyle\prod_{j\in\N^*} \left( \Z/\Z p^{\round nj} \right)^B \\
  \displaystyle\prod_{\on{j\in \N^*}{\beta\in B}} 
    (1+\beta\pi^j)^{m_{j\beta}} U^{(n)} & \mapsto 
  & \left( 
      m_{j\beta} + \Z p^{\round nj} 
    \right)_{\on{j\in \N^*}{\beta\in B}}. 
\end{map}
Passing to projective limits on both sides leads us to the isomorphism of free $\Z_p$-modules
\begin{map}{\mu} 
  U^{(1)} & \to & \M := \Z_p^{\N^*\times B},
\end{map} 
where as usual $\Z_p\subseteq\Q_p$ denotes the ring of $p$-adic integers (\cf Hasse's One-Unit Theorem in \cite[p. 227]{H}). 
Let ${}^{(n)} : \M \onto \M^{(n)}$ be the canonical projection, then $\mu$ can be characterized as the unique map from $U^{(1)}$ to $\M$ satisfying $\mu(y)^{(n)} = \mu^{(n)}(yU^{(n)})$ for all 
$y\in U^{(1)}$ and all $n$. 

For $m = (m_{j\beta})_{\on{j\in \N^*}{\beta\in B}} \in 
\M\setminus\{0\}$ we put 
\[
  \nu(m) := \min \{ 
    jp^{v_p(m_{j\beta})} \sothat j\in\N^*, \beta\in B 
  \}.
\]
Because for $j\in\N^*$ and $\beta \in B$ we have 
$m_{j\beta}^{(n)}=0$ iff $jp^{v_p(m_{j\beta})} \geq n$, 
we can read $\nu(m)$ from $m^{(n)}$ as soon as $n>\nu(m)$.
Let $i := \nu(m)^*$ be the non-$p$-part of $\nu(m)$, then there exists 
$\alpha \in B$ such that 
$ip^{v_p(m_{i\alpha})} = \nu(m) \leq jp^{v_p(m_{j\beta})}$
for all $j \in \N^*$ and $\beta \in B$.
Using~\eqref{roundineq} we see that for each $n$ the logarithmic order of the cyclic group generated by $m^{(n)}$ is 
\begin{eqnarray*}
  \log_p \left| \Z m^{(n)} \right|
  &=& \max \{ 0, \round n{j} - v_p(m_{j\beta}) 
            \sothat j \in \N^*, \beta \in B \} \\
  &=& \max \{ 0, \round n{i} - v_p(m_{i\alpha}) \} = \round n{\nu(m)}.
\end{eqnarray*}

Now let $r\in\N_0$ and $m_i = 
(m_{ij\beta})_{\on{j\in \N^*}{\beta\in B}} \in \M$ 
for $1 \leq i \leq r$ be given.
In the following we describe a method to determine the order of the 
subgroup $\sum_{i=1}^r \Z m_i^{(n)} \subseteq \M^{(n)}$ for all 
$n\in\N$ simultaneously. 
We can assume \wl that $\thru m1r$ are $\Z_p$-linearly independent 
and that $n_1 := \min \{ \nu(m_i) \sothat 1 \leq i \leq r \} = \nu(m_1)$.
As above, for $j_1 := n_1^* \in \N^*$ there is $\beta_1 \in B$ such 
that $v_p(m_{1j_1\beta_1}) = v_p(n_1) \leq v_p(m_{ij_1\beta_1})$ for all 
$i\in\set1r$. 
Replacing $m_i$ by 
$\ti m_i := m_i - \frac{m_{ij_1\beta_1}}{m_{1j_1\beta_1}} m_1$ 
for $i = 2,\ldots,r$ we obtain
\[
  \sum_{i=1}^r \Z m_i^{(n)}
  = \Z m_1^{(n)} \oplus \sum_{i=2}^r \Z \ti m_i^{(n)} 
\]
and $\log_p |\Z m_1^{(n)}| = \round n{n_1}$. 
We apply the same procedure to $\thru{\ti m}2r$ and so on.
In this way we obtain $n_1 \leq\ldots\leq n_r$ such that
\begin{equation} \label{ordereq}
  \log_p \left| \sum_{i=1}^r \Z m_i^{(n)} \right|
  = \sum_{i=1}^r \round n{n_i}
  \quad \forall n\in\N.
\end{equation}
As is evident from the construction, in order to determine $\thru n1r$ 
it suffices to perform the necessary linear transformations on $\thru 
{m^{(n)}}1r$ for an arbitrary $n>n_r$ rather than on $\thru m1r$.

Now suppose that we have found a basis $\thru y1r$ for the group $U_S$ 
with $r=|S|-1\in\N_0$. 
Then by \cite{Ki} we know that $\thru y1r$ are also $\Z_p$-linearly 
independent. 
Apply the described procedure to $m_1:=\mu(y_1),\ldots,m_r:=\mu(y_r)$ to 
obtain $\thru n1r$ as above and define the polynomial 
\[
  \delta_S := \sum_{i=1}^r t^{n_i} \in \Z[t].
\]
Combining equations~\eqref{rounditerate} and~\eqref{ordereq} with the exact sequence~\eqref{useq} results in the 
following recursive formula for the logarithmic orders $\la_S^{(n)} = 
\log_p(U^{(1)}:U_SU^{(n)})$, which at the same time shows that $\delta_S$ 
depends on $\p$ and $S$ only. 
\begin{thm} \label{degthm}
  Write $\delta_S = \sum d_n t^n$ and let $\delta_S^{(n)} := 
  \sum_{l=0}^{v_p(n)} d_{n/p^l} = |\{ i \sothat 1 \leq i \leq r, n_i 
  \leq n, n_i^*=n^* \}|$. 
  Then 
  \[
    \la_S^{(n+1)} = \la_S^{(n)} + e - \delta_S^{(n)}
  \]
  for all $n\in\N$.
\end{thm}
We call $\delta_S$ the \df{$S$-description (at $\p$)} since together with $h_S$ it carries the complete information on the degrees and
the genera of {\em all} ray class field extensions $K_S^{n\p}|K$. 
While the mentioned `first construction' in \cite{XN} was restricted to the special case $q=p$ and $\delta_S = \sum_j t^j$ with the sum over the first $r$ elements of $\N^*$, we are now in the position to deal with practically any choice of $S$ and $\p$. 
\begin{ex}
  We continue Example~\ref{ex020205} and retain the notation therein.
  As we saw, the first $r \in \{1,2,3\}$ of the functions 
  $y_1 = x^{13}(x+1)^6 y^{-3}(y+x^2)^{-2}$,
  $y_2 = x^5 (x+1)^3 y^{-1}(y+x^2)^{-1}$ and 
  $y_3 = x^{-5}(x+1)^{-2}y(y+x^2)$ generate $U_{S_r}$.
  Taking $\p=\p_0$ and $\pi := x$ as the uniformizer at $\p$ yields 
  the expansion $y = \pi^3 + \pi^5 + \pi^6 + \pi^{10} + \cdots$, 
  hence  $\mu^{(6)}(y/\pi^3) = 
  (2,1,1)$ and $\mu^{(6)}((y+x^2)/\pi^2) = (1,1,0)$, 
  where we write $(m_1,m_3,m_5)$ for short instead of 
  $(m_1+8\Z,m_3+2\Z,m_5+2\Z) \in \M^{(6)}$. 
  We obtain the matrix
  \[
    \begin{array}({c}) 
      \mu^{(6)}(y_1) \\ \mu^{(6)}(y_2) \\ \mu^{(6)}(y_3)
    \end{array} = 
    \begin{array}({ccc})
      6&1&1\\
      0&0&1\\
      1&0&1
    \end{array},
  \] 
  from which we read the descriptions $\delta_{S_2} 
  = t^2 + t^5$ and $\delta_{S_3} = t + t^3 + t^5$. 
  Theorem~\ref{degthm} yields 
  $(\la_{S_2}^{(n)})_{n\in\N} = (0,1,1,2,2,2,3,4,4,5,5,6,7,8,\ldots)$ 
  and $(\la_{S_3}^{(n)})_{n\in\N} = (0,0,0,0,0,0,0,1,1,2,\ldots)$, 
  so we see that $\p_3$ does not split completely in $K_{S_2}^{n\p}$ 
  for $n\geq2$, and formulas~\eqref{genfor} and~\eqref{numfor} lead 
  to the following table. 
  \\[\tablesep]
  \begin{flattable*}{10}
         $l$       & 1 &  2 &  3 &  4 &   5 &   6 &   7 &   8 &   9 \\ 
    \hline
    $g(L_{l,S_2})$ & 4 & 10 & 28 & 68 & 164 & 388 & 868 &1892 &4068 \\ 
    \hline
    $N(L_{l,S_2})$ & 7 & 13 & 25 & 49 &  97 & 193 & 385 & 769 &1537 \\ 
  \end{flattable*}
  \\[\tablesep]
  By comparison with the Oesterl\'e bound we see that $L_{2,S_2}$ is 
  optimal. Moreover compared to previous knowledge we have the 
  improved estimates $25 \leq N_2(28) \leq 26$ and 
  $49 \leq N_2(68) \leq 51$. 
\end{ex}
In the same manner we can calculate numerous examples. 
The essence of these is exposed in the tables below. 

Of course we can apply Method~(B) to rational function fields as well, 
which is recommended where Method~(A) fails. 
(Recall that Theorem~\ref{ratthm} only ensures $\delta_S = \sum_n e_S^{(n)} t^n$ 
provided that $pn_S' \geq n_S$.)
\Eg for $q=16$ in addition to 25 possible descriptions which can be 
established by Theorem~\ref{ratthm} we obtain 12 more using Method~(B), namely 
$2t + t^2$, 
$2t + t^3$,
$3t + t^2$,
$3t + t^3$,
$3t + 2t^3$,
$3t + t^2 + t^3$,
$3t + t^2 + 2t^3$,
$3t + t^2 + 2t^3 + t^7$,
$4t + 3t^3 + t^7$,
$4t + 3t^3 + t^5 + t^7$,
$4t + 3t^3 + 2t^5 + t^7$ and
$4t + 3t^3 + 2t^5 + 2t^7$.

%% file: tables.tex
\subsection*{Tables}

Method~(B) and the linear algebra necessary to compute the $e_S^{(n)}$ of 
Method (A) have been implemented in KASH, the KANT shell programming 
language (see \cite{K} for details). 
I am very thankful to the KANT group, and especially to F. He\ss, for 
supporting me in my first steps with their extremely useful number theory 
tool.

The tables below give the currently known range for the number $N_q(g)$ 
for some values of $q$ and $g$ in the form {\em lower bound} -- {\em upper 
bound}, or a single value if $N_q(g)$ is known exactly. 
As a basis we use the results and tables of \cite{GV7}, \cite{Lt3}, 
\cite{Ni}, \cite{NX4} \cite{NX8}, \cite{NX9}, \cite{NX10}, \cite{NX11}, 
\cite{NX12}, \cite{S3}, \cite{XN} and the estimates discussed in Section~2. 

Each entered lower bound is obtained by a corresponding field $L_{l,S}$ for 
an appropriate choice of $K$, $\p$, $l$ and a set $S$ consisting of 
rational places only. 
We bold-face this lower bound where its previously known value is improved 
by our example. 
The entry for $n$ is chosen minimal with $l\leq\la_S^{(n)}$, \ie 
$n = f(L_{l,S},\p)$ is the conductor exponent. 
Also note that we always have $N(L_{l,S}) = h_S(p^l|S|+\eps)$ with $\eps=0$ 
or 1, which is guaranteed by explicitly verifying the supplementary 
condition to inequality~\eqref{numfor}. 

All examples for $q\geq5$ arise from the rational function field and are 
usually computed by means of the quicker Method~(A).
For $q=2$, 3 and 4, where we have also considered ray class fields of other 
than the rational function field, we add two more columns to the table, 
one for the $S$-class number $h_S$ and one for the genus $g(K)$ of the used 
ground field $K$, which in most cases is borrowed from \cite{NX1} and 
\cite{NX3}. 
More detailed tables giving precise reference to how each example using
Method~(B) is obtained are found in~\cite{Au}. 

{\parindent0em
\twocolumn
\begin{figure}
\begin{bigtable}{2}
   6 & 10 & 2 & 1 & 1 & 5 & 1 \\ 
   8 & 11 & 10 & 1 & 5 & 1 & 2 \\ 
  10 & 13 & 4 & 2 & 3 & 1 & 2 \\ 
  13 & 15 & 12 & 1 & 7 & 1 & 4 \\ 
  14 & 15--16 & 14 & 1 & 7 & 1 & 4 \\ 
  15 & 17 & 7 & 3 & 2 & 1 & 0 \\ 
  16 & {\bf17}--18 & 14 & 1 & 8 & 1 & 5 \\ 
  17 & 17--18 & 5 & 4 & 1 & 1 & 0 \\ 
  20 & 19--21 & 0 & 0 & 1 & 19 & 2 \\ 
  22 & 21--22 & 12 & 2 & 5 & 1 & 2 \\ 
  28 & {\bf25}--26 & 7 & 3 & 3 & 1 & 2 \\ 
  29 & 25--27 & 14 & 2 & 6 & 1 & 3 \\ 
  30 & {\bf25}--27 & 12 & 2 & 6 & 1 & 4 \\
  35 & {\bf29}--31 & 16 & 2 & 7 & 1 & 4 \\
  37 & {\bf29}--32 & 14 & 2 & 7 & 1 & 5 \\
  39 & 33 & 8 & 4 & 2 & 1 & 0 \\ 
  41 & {\bf33}--35 & 6 & 3 & 4 & 1 & 4 \\
  42 & {\bf33}--35 & 18 & 2 & 8 & 1 & 5 \\
  44 & {\bf33}--37 & 11 & 3 & 4 & 1 & 2 \\
  55 & {\bf41}--43 & 12 & 3 & 5 & 1 & 3 \\
  58 & {\bf41}--45 & 12 & 3 & 5 & 1 & 4 \\
  60 & {\bf41}--47 & 11 & 3 & 5 & 1 & 4 \\
  61 & {\bf41}--47 & 14 & 3 & 5 & 1 & 3 \\
  68 & {\bf49}--51 & 8 & 4 & 3 & 1 & 2 \\ 
  69 & 49--52 & 16 & 3 & 6 & 1 & 3 \\ 
  71 & {\bf49}--53 & 6 & 4 & 3 & 1 & 3 \\
  74 & {\bf49}--55 & 16 & 3 & 6 & 1 & 4 \\
  75 & {\bf49}--56 & 14 & 3 & 6 & 1 & 5 \\
  79 & 52--58 & 0 & 0 & 2 & 26 & 4 \\
  81 & {\bf49}--59 & 15 & 3 & 6 & 1 & 5 \\
  83 & {\bf57}--60 & 18 & 3 & 7 & 1 & 4 \\
  84 & 57--61 & 18 & 3 & 7 & 1 & 4 \\ 
  89 & {\bf57}--64 & 18 & 3 & 7 & 1 & 5 \\
  91 & {\bf57}--65 & 18 & 3 & 7 & 1 & 5 \\
  95 & 65--68 & 14 & 4 & 4 & 1 & 1 \\ 
  97 & {\bf65}--69 & 8 & 4 & 4 & 1 & 4 \\
  98 & {\bf65}--69 & 20 & 3 & 8 & 1 & 5 \\
 100 & {\bf65}--70 & 12 & 4 & 4 & 1 & 2 \\
\end{bigtable}
\caption{Table for $N_2(g)$}
\end{figure}

\begin{figure}
\begin{bigtable}{3}
  5 & 12--13 & 0 & 0 & 6 & 2 & 3 \\
  7 & 16--17 & 3 & 1 & 5 & 1 & 2 \\
  9 & 19 & 5 & 1 & 6 & 1 & 2 \\
  10 & 19--21 & 5 & 2 & 2 & 1 & 0 \\ 
  13 & 24--25 & 0 & 0 & 2 & 12 & 2 \\ 
  15 & 28 & 6 & 2 & 3 & 1 & 0 \\ 
  17 & 24--30 & 5 & 1 & 4 & 2 & 2 \\
  19 & 28--32 & 12 & 1 & 9 & 1 & 3 \\
  24 & {\bf31}--38 & 14 & 1 & 10 & 1 & 4 \\
  30 & 37--46 & 8 & 2 & 4 & 1 & 1 \\
  33 & {\bf46}--49 & 4 & 2 & 5 & 1 & 3 \\
  36 & 46--52 & 9 & 2 & 5 & 1 & 1 \\ 
  43 & 55--60 & 11 & 2 & 6 & 1 & 1 \\ 
  46 & 55--63 & 6 & 3 & 2 & 1 & 0 \\ 
  47 & {\bf54}--65 & 6 & 2 & 3 & 2 & 1 \\ 
  48 & 55--66 & 11 & 2 & 6 & 1 & 2 \\
  52 & {\bf55}--70 & 11 & 2 & 6 & 1 & 2 \\
  54 & {\bf55}--72 & 9 & 2 & 6 & 1 & 3 \\
  55 & {\bf64}--73 & 12 & 2 & 7 & 1 & 2 \\ 
  57 & {\bf70}--75 & 2 & 1 & 3 & 7 & 3 \\ 
  58 & {\bf55}--77 & 11 & 2 & 6 & 1 & 3 \\
  60 & {\bf64}--79 & 11 & 2 & 7 & 1 & 3 \\
  63 & {\bf64}--82 & 12 & 2 & 7 & 1 & 3 \\
  64 & {\bf64}--83 & 12 & 2 & 7 & 1 & 3 \\
  65 & {\bf72}--84 & 6 & 2 & 4 & 2 & 2 \\
  66 & {\bf64}--85 & 12 & 2 & 7 & 1 & 3 \\
  69 & 82--88 & 8 & 3 & 3 & 1 & 0 \\ 
  70 & {\bf64}--89 & 11 & 2 & 7 & 1 & 4 \\
  72 & {\bf64}--91 & 12 & 2 & 7 & 1 & 4 \\
  73 & {\bf82}--92 & 14 & 2 & 9 & 1 & 3 \\
  75 & {\bf73}--95 & 12 & 2 & 8 & 1 & 4 \\
  78 & {\bf73}--98 & 14 & 2 & 8 & 1 & 4 \\
  79 & {\bf73}--99 & 14 & 2 & 8 & 1 & 4 \\
  81 & {\bf82}--101 & 14 & 2 & 9 & 1 & 4 \\
  82 & {\bf73}--102 & 14 & 2 & 8 & 1 & 4 \\
  84 & {\bf82}--104 & 15 & 2 & 9 & 1 & 4 \\
  85 & {\bf84}--105 & 6 & 1 & 4 & 7 & 3 \\
  87 & {\bf82}--107 & 5 & 3 & 3 & 1 & 2 \\
  93 & {\bf91}--113 & 17 & 2 & 10 & 1 & 4 \\
  94 & {\bf91}--114 & 17 & 2 & 10 & 1 & 4 \\
  96 & {\bf82}--116 & 8 & 3 & 3 & 1 & 1 \\ 
  99 & {\bf100}--119 & 18 & 2 & 11 & 1 & 4 \\ 
\end{bigtable}
\caption{Table for $N_3(g)$}
\end{figure}
\onecolumn
\begin{figure}
\begin{bigtable}{4}
  3 & 14 & 2 & 1 & 3 & 2 & 1 \\
  4 & 15 & 6 & 1 & 7 & 1 & 1 \\ 
  5 & 17--18 & 6 & 2 & 4 & 1 & 0 \\ 
  6 & 20 & 0 & 0 & 4 & 5 & 2 \\
  8 & 21--24 & 6 & 1 & 10 & 1 & 3 \\
  9 & 26 & 3 & 2 & 3 & 2 & 1 \\
  10 & 27--28 & 2 & 2 & 2 & 3 & 1 \\ 
  12 & {\bf29}--31 & 8 & 2 & 7 & 1 & 1 \\ 
  13 & 33 & 6 & 3 & 4 & 1 & 0 \\ 
  19 & {\bf37}--43 & 10 & 2 & 9 & 1 & 2 \\ 
  20 & {\bf37}--45 & 10 & 2 & 9 & 1 & 2 \\
  21 & 41--47 & 8 & 2 & 10 & 1 & 3 \\
  22 & {\bf41}--48 & 10 & 2 & 10 & 1 & 3 \\
  23 & {\bf41}--50 & 7 & 3 & 5 & 1 & 1 \\
  25 & 51--53 & 3 & 2 & 4 & 3 & 2 \\ 
  27 & 49--56 & 6 & 4 & 3 & 1 & 0 \\ 
  28 & {\bf51}--58 & 3 & 3 & 2 & 3 & 1 \\ 
  32 & {\bf57}--65 & 10 & 3 & 7 & 1 & 1 \\ 
  33 & 65--66 & 7 & 4 & 4 & 1 & 0 \\ 
  34 & 57--68 & 8 & 3 & 7 & 1 & 2 \\
  41 & 65--78 & 10 & 3 & 8 & 1 & 2 \\
  47 & {\bf73}--87 & 12 & 3 & 9 & 1 & 2 \\ 
  52 & {\bf81}--95 & 11 & 3 & 10 & 1 & 3 \\
  53 & {\bf81}--96 & 12 & 3 & 10 & 1 & 3 \\
  55 & {\bf81}--99 & 8 & 4 & 5 & 1 & 1 \\
  56 & {\bf81}--101 & 12 & 3 & 10 & 1 & 3 \\
  57 & {\bf81}--102 & 6 & 4 & 5 & 1 & 2 \\
  58 & {\bf97}--103 & 8 & 4 & 6 & 1 & 1 \\
\end{bigtable}\begin{bigtable}4
  61 & 99--108 & 4 & 3 & 4 & 3 & 2 \\ 
  62 & {\bf97}--109 & 7 & 4 & 6 & 1 & 2 \\
  63 & {\bf89}--111 & 14 & 3 & 11 & 1 & 3 \\
  64 & {\bf99}--112 & 3 & 4 & 2 & 3 & 1 \\ 
  65 & 98--114 & 5 & 4 & 3 & 2 & 1 \\
  66 & {\bf97}--115 & 10 & 4 & 6 & 1 & 1 \\
  67 & {\bf97}--117 & 7 & 5 & 3 & 1 & 0 \\ 
  68 & {\bf97}--118 & 10 & 4 & 6 & 1 & 1 \\
  71 & {\bf97}--123 & 8 & 4 & 6 & 1 & 2 \\
  72 & {\bf113}--124 & 10 & 4 & 7 & 1 & 1 \\ 
  73 & {\bf113}--125 & 10 & 4 & 7 & 1 & 1 \\
  74 & {\bf97}--127 & 8 & 4 & 6 & 1 & 2 \\
  76 & 99--130 & 3 & 3 & 4 & 3 & 3 \\ 
  78 & {\bf97}--133 & 7 & 4 & 6 & 1 & 3 \\
  80 & {\bf113}--135 & 10 & 4 & 7 & 1 & 2 \\
  81 & 129--137 & 8 & 5 & 4 & 1 & 0 \\ 
  82 & {\bf113}--138 & 10 & 4 & 7 & 1 & 2 \\
  83 & {\bf113}--140 & 10 & 4 & 7 & 1 & 2 \\
  84 & {\bf113}--141 & 10 & 4 & 7 & 1 & 2 \\
  86 & {\bf113}--144 & 10 & 4 & 7 & 1 & 2 \\
  89 & {\bf113}--148 & 10 & 4 & 7 & 1 & 2 \\
  94 & 129--155 & 11 & 4 & 8 & 1 & 2 \\
  95 & {\bf113}--156 & 10 & 4 & 7 & 1 & 3 \\
  97 & {\bf129}--159 & 12 & 4 & 8 & 1 & 2 \\
  98 & {\bf129}--161 & 7 & 5 & 4 & 1 & 1 \\ 
  99 & {\bf129}--162 & 12 & 4 & 8 & 1 & 2 \\
 100 & {\bf129}--163 & 10 & 4 & 8 & 1 & 3 \\
\end{bigtable}
\caption{Table for $N_4(g)$}
\end{figure}

\begin{figure}\mbox{
\begin{slimtable}5
    6 &       21--25 & 5 & 1 & 4 \\
   22 &       51--60 & 4 & 2 & 2 \\
   34 &       76--83 & 5 & 2 & 3 \\
   56 &     101--125 & 7 & 2 & 4 \\
   70 &     126--150 & 8 & 2 & 5 \\
\mergetable7
    6 &  {\bf22}--32 & 4 & 1 & 3 \\
    9 &  {\bf29}--42 & 5 & 1 & 4 \\
   12 &  {\bf36}--52 & 6 & 1 & 5 \\
   15 &  {\bf43}--60 & 7 & 1 & 6 \\
   45 & {\bf99}--140 & 4 & 2 & 2 \\
   69 &{\bf148}--198 & 5 & 2 & 3 \\
   93 &{\bf197}--254 & 6 & 2 & 4 \\
\end{slimtable}\hspace{0.5em}\begin{slimtable}8
    6 &   33--36 & 6 & 2 & 8 \\
    7 &   33--39 & 4 & 3 & 4 \\
   14 &       65 & 6 & 3 & 8 \\
   29 &  97--109 & 6 & 4 & 6 \\
   38 & 129--135 & 8 & 4 & 8 \\
   46 & 129--158 & 6 & 5 & 4 \\
\mergetable9
    1 &       16 & 3 & 1 & 5 \\
    3 &       28 & 5 & 1 & 9 \\
   12 &   55--63 & 5 & 2 & 6 \\
   15 &   64--74 & 6 & 2 & 7 \\
   21 &   82--95 & 8 & 2 & 9 \\
   48 & 163--180 & 6 & 3 & 6 \\
\end{slimtable}\hspace{0.5em}\begin{slimtable}{16}
    1 &            25 & 4 & 1 & 12 \\
    2 &            33 & 6 & 1 & 16 \\
    5 &        49--55 & 6 & 2 & 12 \\
    6 &            65 & 6 & 2 & 16 \\
   13 &       97--103 & 6 & 3 & 12 \\
   14 &       97--108 & 6 & 3 & 12 \\
   18 &      113--129 & 8 & 3 & 14 \\
   22 &      129--150 &10 & 3 & 16 \\
   27 & {\bf145}--176 & 6 & 4 &  9 \\
   29 &      161--187 & 6 & 4 & 10 \\
   30 &      161--192 & 6 & 4 & 10 \\
   33 & {\bf193}--207 & 7 & 4 & 12\\
   38 & {\bf193}--233 & 8 & 4 & 12 \\
   42 & {\bf209}--254 & 8 & 4 & 13 \\
   50 & {\bf225}--291 &10 & 4 & 14 \\
\end{slimtable}
}\caption{Tables for $N_q(g)$ with $q=5$, 7, 8, 9 and 16}
\end{figure}

\begin{figure}\mbox{
\begin{slimtable}{25}
    2 &            46 & 3 & 1 &  9 \\
    4 &            66 & 4 & 1 & 13 \\
    6 &   {\bf66}--85 & 5 & 1 & 13 \\
   10 &           126 & 7 & 1 & 25 \\
   12 & {\bf101}--140 & 3 & 2 &  4 \\
   22 & {\bf176}--208 & 4 & 2 &  7 \\
   24 & {\bf176}--222 & 4 & 2 &  7 \\
   32 & {\bf226}--275 & 5 & 2 &  9 \\
   34 & {\bf226}--289 & 5 & 2 &  9 \\
   36 & {\bf226}--302 & 5 & 2 &  9 \\
\mergetable{27}
    4 &        64--68 & 3 & 2 &  7 \\
   10 &  {\bf91}--128 & 5 & 2 & 10 \\
   12 &      109--148 & 5 & 2 & 12 \\
   13 &      136--156 & 3 & 3 &  5 \\
   15 & {\bf136}--171 & 6 & 2 & 15 \\
   16 & {\bf136}--178 & 6 & 2 & 15 \\
   21 &      163--214 & 8 & 2 & 18 \\
   36 &      244--319 &11 & 2 & 27 \\
   39 &      271--340 & 5 & 3 & 10 \\
   48 & {\bf325}--402 & 6 & 3 & 12 \\
\end{slimtable}\hspace{0.3em}\begin{slimtable}{32} 
    6 &        81--99 &  6 & 2 & 20 \\
   11 & {\bf113}--154 &  6 & 3 & 14 \\
   12 &      129--165 & 10 & 2 & 32 \\
   13 & {\bf129}--176 &  6 & 3 & 16 \\
   14 &      145--187 &  6 & 3 & 18 \\
   18 & {\bf161}--220 &  8 & 3 & 20 \\
   27 & {\bf209}--290 & 10 & 3 & 26 \\
   28 &      257--298 & 10 & 3 & 32 \\
   29 & {\bf225}--306 &  6 & 4 & 14 \\
   30 &      273--313 &  6 & 4 & 17 \\
   38 & {\bf289}--375 &  8 & 4 & 18 \\
\mergetable{81}
   1 &           100 &  3 & 1 & 33 \\
   3 &           136 &  5 & 1 & 45 \\
   4 &           154 &  6 & 1 & 51 \\
   9 &           244 & 11 & 1 & 81 \\
  12 &           298 &  5 & 2 & 33 \\
  15 & {\bf280}--352 &  6 & 2 & 31 \\
  21 & {\bf352}--459 &  8 & 2 & 39 \\
  31 & {\bf460}--639 & 11 & 2 & 51 \\
  36 &           730 & 11 & 2 & 81 \\
  48 & {\bf676}--885 &  6 & 3 & 25 \\
\end{slimtable}\hspace{0.3em}\begin{slimtable}{49}
    3 &            92 &  3 & 1 & 13 \\ 
    9 &           176 &  5 & 1 & 25 \\
   21 &           344 &  9 & 1 & 49 \\
   45 & {\bf442}--579 &  4 & 2 &  9 \\
\mergetable{64}
    1 &            81 &  4 & 1 & 40 \\
    2 &            97 &  6 & 1 & 48 \\
    4 &           129 & 10 & 1 & 64 \\
    6 &           161 &  6 & 2 & 40 \\
    8 & {\bf169}--193 &  8 & 2 & 42 \\
   11 &      201--241 & 10 & 2 & 50 \\
   12 &           257 & 10 & 2 & 64 \\
   13 & {\bf225}--272 &  6 & 3 & 28 \\
   14 & {\bf241}--288 &  6 & 3 & 30 \\
   18 & {\bf281}--352 &  8 & 3 & 35 \\
   22 & {\bf321}--416 & 10 & 3 & 40 \\
   24 & {\bf337}--448 & 10 & 3 & 42 \\
   26 & {\bf385}--480 & 10 & 3 & 48 \\
   27 &      401--496 & 10 & 3 & 50 \\
   28 &           513 & 10 & 3 & 64 \\
   30 & {\bf385}--536 &  6 & 4 & 24 \\
   38 & {\bf449}--627 &  8 & 4 & 28 \\
   50 & {\bf561}--762 & 10 & 4 & 35 \\
\end{slimtable}
}\caption{Tables for $N_q(g)$ with $q=25$, 27, 32, 49, 64 and 81}
\end{figure}
}

%% file: rami.tex
\section*{Appendix: Ramification Groups}

Let $L|K$ be a finite Galois extension of global (or local) fields with 
group $G(L|K)$, $\p$ a (non-archimedian) place of $K$ and $\q$ a place of 
$L$ lying above $\p$. 
For $s\in[-1,\infty)$ we have the \df{lower ramification groups of 
$\q|\p$}, 
\[
  G_s(\q|\p) := \{ \sigma\in G(L|K) \sothat v_\q(\sigma z-z)\geq s+1
                   \mbox{ for all $z\in L$ integral at }\q \},
\]
which are renumbered to define the \df{upper ramification groups} 
$G^{\eta(s)}(\q|\p) := G_s(\q|\p)$ by means of the continuous bijection
\begin{map}{\eta}
  [-1,\infty) & \to & [-1,\infty) \\
        s & \mapsto & \displaystyle \int_0^s 
                      \frac{|G_x(\q|\p)|}{|G_0(\q|\p)|} dx
\end{map}
(\cf \cite[pp.\ 186ff]{Nk}).
Note that $\eta$ is the identity at least on $[-1,0]$ and that 
$G_{-1}(\q|\p) = G^{-1}(\q|\p)$ is the decomposition group and 
$G_0(\q|\p) = G^0(\q|\p)$ is the inertia group of $\q|\p$. 

Now let $L|K$ be abelian. 
Then the ramification groups do no longer depend on the specific choice of  
$\q$, so we can write $G^s(L,\p)$ instead of $G^s(\q|\p)$. 
In Section 1 we made use of an upper index version of Hilbert's 
\renewcommand\thethm{}
\begin{diffor}
  Denote by $d(\q|\p)$ the different exponent of $\q|\p$. It satisfies
  \[
    d(\q|\p) = \sum_{n=0}^\infty \left( 
                 |G^0(L,\p)| - (G^0(L,\p):G^n(L,\p)) 
               \right)
  \]
\end{diffor}

\begin{proof}
  For convenience put $G^s:=G^s(L,\p)$ and $\psi:=\eta^{-1}$, where 
  $\eta$ is defined as above. 
  Let $t_0:=-1$, and let $0\leq t_1 < \ldots < t_r$ be the (other) 
  \df{jumps} of the filtration $(G^s)_{s\geq-1}$, \ie the numbers $s\geq-1$ 
  satisfying $G^{s+\eps} \subsetneq G^s$ for all $\eps>0$ (\cf \cite[p.\ 
  372]{Nk}). 
  Obviously all $\psi(t_i)$ are integers, and $t_0=-1$ is a jump iff the 
  inertia degree of $\q|\p$ is $>1$. The bijection $\eta$ is piecewise 
  linear, its slopes being 
  \[
    \eta'(s) = \frac1{(G^0:G^{t_i})} 
             = \frac{t_i-t_{i-1}}{\psi(t_i)-\psi(t_{i-1})}
  \]
  for $\psi(t_{i-1})<s<\psi(t_i)$ and $i\in\set1r$.
  By the Hasse-Arf theorem the $t_i$ are integers, too. 
  Substituting into Hilbert's (lower index) different formula (see 
  \cite[p.\ 124]{St}) yields 
  \begin{eqnarray*}
    d(\q|\p) &=& \sum_{n=0}^\infty (|G_n(\q|\p)|-1) 
              =  \sum_{i=1}^r \sum_{s=\psi(t_{i-1})+1}^{\psi(t_i)} 
                   (|G^{t_i}|-1) \\
             &=& \sum_{i=1}^r (\psi(t_i) - \psi(t_{i-1}))
                   \frac {|G^0| - (G^0:G^{t_i})} {(G^0:G^{t_i})} \\
             &=& \sum_{i=1}^r (t_i - t_{i-1}) (|G^0| - (G^0:G^{t_i})) \\
             &=& \sum_{i=1}^r \sum_{n=t_{i-1}+1}^{t_i}                
                   (|G^0| - (G^0:G^n))
              =  \sum_{n=0}^\infty (|G^0| - (G^0:G^n)). \qed
  \end{eqnarray*}  
  \noqed
\end{proof}

%% file: address.tex
Afd.\ wiskunde \\
Rijksuniversiteit Groningen \\
Blauwborgje 3 \\ 
NL-9747 AC Groningen \\
e-mail: auer@math.rug.nl